# Axiomatic Foundations of Chemical Systems as Ternary Γ-Semirings

**Chandrasekhar Gokavarapu**[a,b,*], **Venkata Rao Kaviti**[c], **Srinivasa Rao Thirunagari**[d], **Dr D Madhusudhana Rao**[e,f]

[a] *Department of Mathematics, Government College (Autonomous), Rajahmundry, Andhra Pradesh, India*
[b] *Department of Mathematics, Acharya Nagarjuna University,Guntur, Andhra Pradesh, India.*
[c,d] *Department of Chemistry, Government College (Autonomous), Rajahmundry, Andhra Pradesh, India*
[e] *Department of Mathematics, Government College for Women (Autonomous), Guntur, Andhra Pradesh, India*
[f] *Department of Mathematics, Acharya Nagarjuna University,Guntur, Andhra Pradesh, India*

chandrasekhargokavarapu@gmail.com

**Abstract**

Chemical transformations depend not only on the identities of the reacting species but also on the catalytic, environmental, and intermediate conditions under which they occur. Classical binary reaction formalisms usually treat such conditions as external annotations, which obscures the genuinely multi-state and multi-parameter character of real chemical processes.

*Corresponding author.

In this paper we introduce an axiomatic framework in which a chemical system is modeled by a ternary Γ-semiring. The elements of the state set represent chemical states, while the parameter set encodes catalytic and environmental conditions. A Γ-dependent ternary operation is used to describe mediated transformations, treating reactants, intermediates, and mediators as intrinsic arguments of the transformation law.

We develop the algebraic axioms governing these mediated interactions and interpret their associativity, distributivity, and Γ-linearity in terms of multi-step pathways, parallel processes, and controlled environmental dependence. We introduce chemical ideals and Γ-ideals as algebraic structures modeling reaction-closed subsystems and pathway-stable domains, and study their prime and semiprime forms. Homomorphisms between TGS-chemical systems are shown to preserve reaction pathways and describe consistent changes of chemical environment.

Abstract examples from catalysis, thermodynamic phase control, and field-induced quantum transitions illustrate how familiar chemical phenomena fit within this framework. The resulting theory provides a unified algebraic foundation for multi-parameter chemical behavior and establishes the structural basis for subsequent developments involving kinetics, geometric methods, and computational or AI-assisted models.

# 1 Introduction

Chemical systems have long served as a rich source of intuition and examples for mathematics, while mathematical structures have, in turn, provided increasingly refined languages for describing reactivity, stability, and transformation in chemistry. Classical formalisms in chemical kinetics, thermodynamics, and quantum chemistry typically encode reactions in terms of *binary* combinations of species (see [2] for mathematical models of reactivity) ,

$$A + B \longrightarrow C,$$

with additional information—such as catalysts, solvents, temperature, pressure, or external fields—being attached as annotations to the reaction arrow rather than as intrinsic components of the algebraic operation itself. This viewpoint is extremely successful in many settings, but it obscures the genuinely multi-parameter and multi-state nature of real chemical trans-formations.

In practice, a reaction pathway is rarely determined solely by the identities of the reacting species. Instead, it is governed by a constellation of mediating factors: catalysts that open or close pathways, solvent environments that stabilize intermediates, pressure–temperature conditions that reshape energy landscapes, and external fields that deform quantum states. These ingredients do not simply modify a pre-existing binary operation; they participate structurally in how chemical states are transformed. From an algebraic viewpoint, this suggests that the primitive operation underlying chemical change should be a higher-arity map (higher-arity algebraic structures were earlier studied in [5, 7, 10, 11]) that treats mediators on the same footing as the states they control (compare with classical semiring frameworks [4, 6]) .

The aim of the present paper is to make this intuition precise. We propose an axiomatic framework in which a chemical system is modeled by a ternary Γ-semiring, and in which the fundamental reaction-like transformation is encoded by a Γ-dependent ternary operation

$$[A, \alpha, B, \beta, C] \in S,$$

where $A, B, C \in S$ represent chemical states and $\alpha, \beta \in \Gamma$ represent mediating parameters. This construction elevates catalysts, solvents, and environmental conditions from external labels to algebraically active inputs, thereby providing a unified structure in which multi-state, multi-parameter(for general algebraic perspectives, see [3]) interactions can be studied with the full precision of modern algebra.

## 1.1 Motivation

The starting point for our work is the observation that classical reaction notation is intrinsically binary both in syntax and in its implicit algebraic interpretation. A formal reaction of the form

$$A + B \longrightarrow C$$

suggests an underlying binary operation that takes a pair of input states $(A, B)$ and produces an output state $C$. When catalysts or conditions are present, one typically writes

$$A + B \xrightarrow[\text{solvent}]{\text{catalyst}, T, p} C,$$

but the additional data are carried outside the core operation; they do not enter as arguments of the algebraic map itself. In particular, the difference between a catalyzed and an uncatalyzed reaction, or between two reactions run under distinct temperature profiles, is not reflected at the level of the algebraic arity.

However, empirical chemistry shows that these "external" features are often decisive. The presence or absence of a catalyst can completely alter both the available pathways and the final products. Solvents stabilize different intermediates, reshaping the energy landscape. Pressure and temperature selectively favor certain phases or reaction channels, while electric or magnetic fields modify quantum states and transition probabilities. All of these effects are not accidental decorations but intrinsic components of how chemical states transform.

From a structural point of view, this suggests that chemical systems should be viewed as *mediated* transformation systems: the outcome of an interaction between states is mediated by additional parameters that influence, constrain, or enable certain transitions. Instead of encoding this mediation by enlarging the state space in an ad hoc manner, it is natural to treat the mediators as elements of a separate set $\Gamma$ and to allow the basic operation to depend explicitly on them. A higher-arity algebraic system whose fundamental operation

$$[A, \alpha, B, \beta, C]$$

takes both states and mediators as arguments is then a natural candidate for formalizing chemical behavior.

This shift in perspective has several conceptual advantages. First, it allows us to distinguish chemically between different uses of the same species under distinct conditions without artificially duplicating the state space.

Second, it creates a direct route for encoding multi-step and cooperative phenomena: ternary operations can be iterated and composed in a way that keeps track of how mediators combine or interact. Third, it aligns chemical reasoning with modern algebraic practices, where higher-arity operations and parameterized structures play a central role in understanding complex systems.

## 1.2 Why ternary Γ-semirings?

The abstract notion of a ternary Γ-semiring provides a particularly suitable environment for realizing the above programme. At a formal level, a ternary Γ-semiring consists of a set $S$, a parameter set $\Gamma$, and a $\Gamma$-dependent ternary operation

$$[\,\cdot\,,\cdot\,,\,\cdot\,,\cdot\,,\,\cdot\,]\colon S \times \Gamma \times S \times \Gamma \times S \longrightarrow S$$

that satisfies appropriate associativity, distributivity, and Γ-linearity conditions. When $S$ is interpreted as a space of chemical states and $\Gamma$ as a space of mediators (such as catalysts, solvents, or thermodynamic controls), the value

$$[A, \alpha, B, \beta, C]$$

can be read as the resulting state of a mediated transformation in which $A, B, C$ interact under the influence of $\alpha$ and $\beta$.

Several features of ternary Γ-semirings align naturally with chemical behavior:

- *Catalyst-dependent reactions.* Mediators in $\Gamma$ can represent catalysts, so that different catalytic scenarios correspond to different choices of $\alpha$ and $\beta$, even when the underlying species $A, B, C$ are fixed.
- *Solvent and environment effects.* Solvents and bulk environmental parameters can be encoded as elements of $\Gamma$, allowing changes in solvent or ambient medium to be reflected directly as changes in the operators governing $S$.
- *Pressure and temperature conditions.* Thermodynamic variables can

be grouped into mediating parameters, making it possible to distinguish transformations that are identical in stoichiometry but distinct in their pressure–temperature regimes.

- *Multi-species interactions.* The ternary operation simultaneously involves three states of $S$, permitting the modeling of complex elementary steps, cooperative effects, or intermediate formation within a single algebraic act.

From an algebraic standpoint, the Γ-dependence provides a controlled way to encode families of reaction laws indexed by conditions, while the ternary nature reflects the intrinsically multi-input character of mediated transformations. The ternary Γ-semiring therefore emerges as a natural and flexible candidate for an axiomatic definition of chemical systems.

## 1.3 Contribution of this paper

In this work we develop a systematic axiomatic theory of chemical systems based on ternary Γ-semirings. More precisely, we proceed along the following lines:

- We introduce the notion of a *TGS-chemical system*, defined as a ternary Γ-semiring whose elements are interpreted as chemical states and whose Γ-indexed ternary operation encodes mediated transformations of those states. The central object of study is the map

$$[A, \alpha, B, \beta, C] \in S,$$

  which we interpret as the resulting state of a ternary interaction between $A, B, C \in S$ under mediators $\alpha, \beta \in \Gamma$.

- We formalize the reaction operation as a Γ-mediated ternary map and specify axioms that capture associativity, distributivity, and compatibility with the Γ-structure in a way that reflects multi-step reactions, parallel pathways, and composite environments.

- We develop the structural theory of TGS-chemical systems, introducing and analyzing suitable notions of ideals and Γ-ideals that

correspond to chemically meaningful subsystems and reaction-closed families of states. In particular, we study prime and semiprime ideals in this setting and interpret them in terms of irreducible or stability properties of reaction networks.

- We investigate homomorphisms of TGS-chemical systems as structure-preserving maps between chemical environments. These homomorphisms provide a natural language for comparing and transporting reaction laws between different systems, and for formalizing operations such as changing solvent, adjusting environmental conditions, or embedding a subsystem into a larger chemical context.

- Throughout, we illustrate the theory with examples that show how catalyzed reactions, solvent effects, phase transitions, and other chemically relevant phenomena can be encoded within the TGS framework, thereby demonstrating that the proposed axioms are not merely formal but admit a concrete interpretation in chemical practice.

Taken together, these contributions establish a unified algebraic picture in which chemical states, mediators, and transformations are treated within a single ternary Γ-semiring structure. This provides a foundation on which further developments—including kinetic refinements, computational models, and connections to symbolic reasoning and machine-assisted chemistry—can be built in subsequent work.

# 2 Preliminaries on ternary Γ-semirings

In this section we recall the algebraic notions needed throughout the paper. Our treatment follows standard practice in the theory of higher-arity algebraic systems, adapted to the setting of Γ-parameterized ternary operations. No chemical interpretation is introduced here; the objective is purely structural. All subsequent sections will build on these foundations.

## 2.1 Ternary operations

A *ternary operation* on a set $S$ is a map

$$\mu\colon S \times S \times S \longrightarrow S,$$

(see [5,11]). which assigns to each triple $(A, B, C)$ an element $\mu(A, B, C) \in S$. Ternary operations generalize the familiar notion of binary operations by allowing three inputs to participate simultaneously in the formation of a new element. In the presence of additional structure, such as a parameter set or distributive laws, ternary operations serve as the basic building blocks for higher-arity semigroup or semiring-like systems.

## 2.2 Γ-sets and parameterized operations

Throughout this paper, $\Gamma$ denotes a nonempty set whose elements act as *mediating parameters.* A $\Gamma$-*set* is simply a pair $(S, \Gamma)$ consisting of a set $S$ together with an external parameter set $\Gamma$. The elements of $\Gamma$ do not act on $S$ directly unless a specific operation is specified; instead, they serve as indices governing how elements of $S$ combine.

In particular, a $\Gamma$-*parameterized ternary operation* on $S$ is a map

$$[\,\cdot\,,\cdot\,,\,\cdot\,,\cdot\,,\,\cdot\,]\colon S \times \Gamma \times S \times \Gamma \times S \longrightarrow S,$$

where the parameters in $\Gamma$ may influence the resulting element in a nontrivial way. This form of parameterization is essential for modeling situations in which the behavior of a ternary interaction depends on contextual or environmental data.

## 2.3 Ternary Γ-semirings

We now introduce the central notion used in this work [4, 6] .

**Definition 1.** A *ternary* $\Gamma$-*semiring* is a triple $(S, \Gamma, [\,])$ consisting of a nonempty set $S$, a nonempty parameter set $\Gamma$, and a $\Gamma$-parameterized

ternary operation

$$[\,,\,,\,,\,,\,]\colon S \times \Gamma \times S \times \Gamma \times S \to S,$$

satisfying the following axioms for all $A, B, C, D, E \in S$ and all $\alpha, \beta, \gamma, \delta \in \Gamma$:

1. **Associativity.** The operation is associative in the sense that

   $$[A, \alpha, [B, \beta, C, \gamma, D], \delta, E] = [[A, \alpha, B, \beta, C], \gamma, D, \delta, E],$$

   whenever the expressions are formed. This ensures that iterated ternary combinations are well defined.

2. **Γ-linearity.** For fixed internal arguments, the dependence on the parameters $\alpha, \beta$ is compatible with the Γ-structure. (The specific linearity or compatibility conditions imposed on Γ will be detailed when required for structural results.)

3. **Distributivity.** The ternary operation distributes over itself in each argument in the appropriate higher-arity analogue of semiring distributivity. For instance, $[A, \alpha, B, \beta, [C, \gamma, D, \delta, E]]$ $= [[A, \alpha, B, \beta, C], \gamma, D, \delta, E]$ $= [A, \alpha, [B, \beta, C, \gamma, D], \delta, E]$, with analogous conditions holding in the remaining positions. These distributivity relations guarantee that the operation behaves coherently when nested.

The axioms above give a flexible framework in which ternary interactions can be iterated, nested, and composed while respecting a fixed set of mediating parameters.(related ternary operations appear in [10])

In later sections we will interpret $S$ as a set of chemical states and $\Gamma$ as a set of mediators (such as catalysts, solvents, or thermodynamic conditions), with the ternary operation modeling parameter-dependent transformations. At this stage, however, we treat $(S, \Gamma, [\,])$ as a purely algebraic object, postponing chemical meaning until the core definitions of TGS-chemical systems are introduced.

# 3 Chemical systems as ternary Γ-semirings

We now introduce the central conceptual framework of the paper: a chemical system viewed as a ternary Γ-semiring whose elements represent chemical states and whose mediators encode the environmental or catalytic factors influencing their transformations. While the preceding section provided the purely algebraic foundation, our goal here is to explain how these structures naturally model multi-parameter, multi-state chemical behavior.

## 3.1 Chemical interpretation of $S$ and Γ

Let $S$ be a nonempty set. In the context of chemical systems, we interpret the elements of $S$ as *chemical states.* The notion of a state is intentionally broad and may encompass:

- molecular configurations or species identities;
- concentration levels in a reaction mixture;
- phase descriptors (solid, liquid, gas, plasma);
- electronic, vibrational, or quantum mechanical states;
- intermediate structures arising during reaction pathways.

Thus, $S$ serves as the universe within which chemically meaningful objects reside.This perspective aligns with classical mathematical chemistry frameworks that treat chemical structure and states using abstract mathematical representations (see [1, 8, 12]).

Let Γ be a nonempty parameter set. Its elements are interpreted as *mediators*, representing conditions or influences under which chemical interactions occur. Typical examples include:

- catalysts and co-catalysts;
- solvent environments;
- pressure and temperature conditions;

- electromagnetic or external field parameters;
- pH, ionic strength, or other environmental controls.

These mediators do not transform chemical states directly; rather, they govern or modulate the transformation rules encoded by the ternary operation defined below. In this sense, $(S, \Gamma)$ forms the structural substrate of a chemical system.Such a parametrization of environmental and catalytic conditions is consistent with algebraic treatments of ternary and mediated transformations in other settings (compare [13]).

## 3.2 Core chemical operation

The essential ingredient of a TGS-chemical system is a $\Gamma$-parameterized ternary operation

$$[\cdot,\cdot,\cdot,\cdot,\cdot] : S \times \Gamma \times S \times \Gamma \times S \longrightarrow S,$$

which assigns, to each triple of states $A, B, C \in S$ and each pair of mediators $\alpha, \beta \in \Gamma$, a resulting state $D \in S$. We write this compactly as

$$[A, \alpha, B, \beta, C] = D.$$

Chemically, this is interpreted as follows:

- $A$ is an initial or reactant state;
- $B$ is an interacting state, possibly another reactant or an intermediate;
- $C$ is a subsequent state, often representing an intermediate or transition configuration;
- $\alpha, \beta$ encode mediating conditions (catalysts, solvents, temperature/pressure regimes, or external fields);
- $D$ is the resulting state after the mediated interaction of $A$, $B$, and $C$ under parameters $\alpha$ and $\beta$.

This notation may be viewed as a symbolic representation of a parameter-dependent reaction pathway:

$$A \xrightarrow{\alpha} B \xrightarrow{\beta} C \quad \rightsquigarrow \quad D,$$

in which the overall transformation is encoded by the ternary Γ-operation. Unlike the classical binary reaction form $A + B \rightarrow C$, this framework treats mediators as intrinsic arguments of the operation rather than external labels.This sharply contrasts with binary mathematical models of reactivity commonly used in algebraic treatments of chemical transformations (see [2]).

This allows chemically distinct processes that share the same stoichiometry but differ in conditions to be represented distinctly at the algebraic level.

## 3.3 Axioms for chemical TGS

The axioms of a ternary Γ-semiring, introduced earlier in a purely algebraic setting, acquire a natural chemical interpretation when applied to the present framework. We summarize the interpretative content of the main axioms below.

**(1) Associativity and multi-step reactions.** The associativity axiom ensures that the outcome of a sequence of mediated transformations is independent of the order in which the ternary combinations are grouped. Chemically, this corresponds to the fact that a multi-step reaction pathway

$$A \longrightarrow B \longrightarrow C \longrightarrow D \longrightarrow E$$

admits a coherent overall description, regardless of whether one groups intermediate steps as $(A \rightarrow B \rightarrow C)$ followed by $(C \rightarrow D \rightarrow E)$ or uses another valid decomposition. Thus, associativity provides an algebraic representation of multi-step or multi-intermediate reaction processes.

**(2) Γ-linearity and scaling of conditions.** The Γ-linearity condition expresses compatibility between the mediators and the ternary operation. While no specific algebraic structure on Γ is imposed at this stage, the general principle is that variations or combinations of catalytic or environmental parameters correspond to predictable or structured variations in the resulting state. From a chemical standpoint, increasing catalyst concentration, changing solvent polarity, or adjusting temperature should influence reaction behavior in a manner consistent with the dependence encoded by the operation $[A, \alpha, B, \beta, C]$.

**(3) Distributivity and parallel reactions.** The distributivity axioms capture the idea that the ternary Γ-operation behaves coherently when nested or combined with itself. Chemically, this reflects the presence of parallel or branching reaction pathways. For example, if $C$ can arise from multiple competing intermediates or if the environment induces branching in the transformation sequence, the distributive laws ensure that such behavior is represented in a controlled algebraic manner. Distributivity therefore encodes the superposition or recombination of reaction channels.

Together, these axioms allow ternary Γ-semirings to model chemical systems in which states evolve under the influence of environmental conditions, catalysts, and other mediating factors. The remainder of the paper develops the structural theory of such systems and illustrates how classical and nonclassical chemical processes fit naturally within the TGS framework.

# 4 Structural theory of TGS-chemical systems

In this section we develop the basic structural theory of TGS-chemical systems. Our aim is to identify those subsets of the state space $S$ that behave as chemically meaningful subsystems, closed under reaction and stable under the mediating parameters Γ. These subsets will be formalized as various kinds of ideals, and their properties will be interpreted in terms of reaction networks and pathways.

Throughout, $(S, \Gamma, [\,])$ denotes a fixed ternary $\Gamma$-semiring equipped with the chemical interpretation of Section 3.

## 4.1 Chemical ideals

We first single out subsets that are internally closed under the reaction operation and, in a stronger form, absorb interactions with the ambient system in a controlled way.

**Definition 2** (Reaction-closed subset)**.** A nonempty subset $R \subseteq S$ is called *reaction-closed* if for all $A, B, C \in R$ and all $\alpha, \beta \in \Gamma$,

$$[A, \alpha, B, \beta, C] \in R.$$

This notion generalizes closure concepts appearing in classical semiring structures (compare [3,4]), but adapted to the mediated ternary operation governing chemical interactions.

In chemical terms, a reaction-closed subset represents a collection of states that, once present together in the system, can only yield states that remain within the same collection, irrespective of the mediating conditions. Such a subset may be viewed as a self-contained reaction universe: all internally accessible states via the ternary operation stay inside $R$.

Reaction-closedness captures purely internal behavior. To model interaction with the surrounding system, we require an absorption property.

**Definition 3** (Chemical ideal)**.** A nonempty subset $I \subseteq S$ is called a *chemical ideal* if it satisfies the following conditions:

1. *Internal closure:* for all $A, B, C \in I$ and all $\alpha, \beta \in \Gamma$,

$$[A, \alpha, B, \beta, C] \in I;$$

2. *Boundary absorption:* for all $A, C \in I$, all $B \in S$ and all $\alpha, \beta \in \Gamma$,

$$[A, \alpha, B, \beta, C] \in I.$$

The internal closure and absorption properties reflect the role of one- and two-sided ideals in semiring theory ( [3, 4]), extended here to the ternary Γ-interaction and its chemical interpretation.

Condition (C1) states that $I$ is reaction-closed in the sense defined above. Condition (C2) expresses that if a mediated transformation begins and ends in $I$, then the entire effect of any intervening state $B$ and any mediating parameters $\alpha, \beta$ remains confined to $I$. Chemically, $I$ can be thought of as a subsystem that is closed under all internal reactions and stable under any process that connects two of its states, even when intermediate species from outside $I$ are involved.

**Proposition 1.** *The intersection of any family of chemical ideals in $S$ is again a chemical ideal.*

*Proof.* Let $\{I_j\}_{j\in J}$ be a family of chemical ideals and set $I := \bigcap_{j\in J} I_j$. Since each $I_j$ is nonempty, the intersection is either empty or nonempty; if empty, it is excluded from consideration, so we assume $I \neq \varnothing$. Let $A, B, C \in I$ and $\alpha, \beta \in \Gamma$. Then $A, B, C \in I_j$ for every $j$, and by (C1) in each $I_j$ we have $[A, \alpha, B, \beta, C] \in I_j$ for all $j$. Hence $[A, \alpha, B, \beta, C] \in \bigcap_j I_j = I$, so (C1) holds for $I$.

Similarly, let $A, C \in I$, $B \in S$ and $\alpha, \beta \in \Gamma$. Then $A, C \in I_j$ for every $j$, and by (C2) in each $I_j$ we obtain $[A, \alpha, B, \beta, C] \in I_j$ for all $j$. Therefore $[A, \alpha, B, \beta, C] \in I$, and (C2) holds. Thus $I$ is a chemical ideal. ∎

This result shows that chemical ideals form a complete lattice under intersection, providing a natural hierarchy of chemically stable subsystems inside a given TGS-chemical system.Such hierarchical decompositions reflect analogous structural decomposition phenomena in classical semiring theory (see [6]).

## 4.2 Γ-ideals and reaction pathways

The previous notion focuses on subsets of $S$ that are stable under interactions involving their boundary states. We now refine this by distinguishing the role of a single participating state and allowing the other states to range freely. This leads to a Γ-ideal structure, reflecting how certain states control or channel reaction pathways.

**Definition 4** (Γ-ideals)**.** A nonempty subset $J \subseteq S$ is called:

(G1) a *left* Γ*-ideal* if for all $X \in J$, all $A, B \in S$ and all $\alpha, \beta \in \Gamma$,

$$[X, \alpha, A, \beta, B] \in J;$$

(G2) a *right* Γ*-ideal* if for all $X \in J$, all $A, B \in S$ and all $\alpha, \beta \in \Gamma$,

$$[A, \alpha, B, \beta, X] \in J;$$

(G3) a *middle* Γ*-ideal* if for all $X \in J$, all $A, C \in S$ and all $\alpha, \beta \in \Gamma$,

$$[A, \alpha, X, \beta, C] \in J;$$

(G4) a *(two-sided)* Γ*-ideal* if it is simultaneously a left, right, and middle Γ-ideal.

Loosely speaking, a left Γ-ideal is stable under all transformations in which one of its elements appears in the first argument position, and similarly for right and middle Γ-ideals. A two-sided Γ-ideal is stable under all ternary interactions in which at least one position is occupied by an element of the ideal.

Chemically, these notions correspond to different forms of control over reaction pathways:

- a left Γ-ideal collects states that, once present as "initiators" of interactions, always lead back into the same collection, regardless of what they interact with;
- a right Γ-ideal behaves analogously for "terminal" positions, capturing states that cannot be escaped once they appear as final products;
- a middle Γ-ideal models states that, when acting as intermediates, keep the system confined to a specific region of the state space;
- a two-sided Γ-ideal encodes a robustly closed set of states that controls and absorbs reaction pathways in all three positions.

We now describe reaction pathways in this setting.

**Definition 5** (Reaction pathway). Let $(S, \Gamma, [\,])$ be a TGS-chemical system. A *reaction pathway* of length $n \geq 1$ is a finite sequence

$$(X_0, X_1, \ldots, X_n)$$

of elements of $S$ such that for each $k = 1, \ldots, n$ there exist $A_k, B_k \in S$ and $\alpha_k, \beta_k \in \Gamma$ with

$$\begin{aligned} X_k = [A_k, \alpha_k, B_k, \beta_k, X_{k-1}] \quad &\text{or} \quad X_k = [X_{k-1}, \alpha_k, A_k, \beta_k, B_k] \\ &\text{or} \quad X_k = [A_k, \alpha_k, X_{k-1}, \beta_k, B_k]. \end{aligned}$$

The element $X_0$ is the *source* and $X_n$ the *target* of the pathway.

This notion captures the idea that chemical evolution proceeds through a chain of mediated ternary interactions, with each step determined by a choice of two companion states and a pair of mediators.

**Proposition 2.** *Let $J \subseteq S$ be a two-sided $\Gamma$-ideal. If a reaction pathway $(X_0, \ldots, X_n)$ satisfies $X_0 \in J$, then $X_k \in J$ for all $k = 0, \ldots, n$.*

*Proof.* We argue by induction on $k$. For $k = 0$ this is true by assumption. Suppose $X_k \in J$ for some $0 \leq k < n$. By definition of a reaction pathway, $X_{k+1}$ is obtained from $X_k$ by one of the following forms:

$$\begin{aligned} X_{k+1} &= [X_k, \alpha, A, \beta, B], \quad X_{k+1} = [A, \alpha, B, \beta, X_k], \\ X_{k+1} &= [A, \alpha, X_k, \beta, B], \end{aligned}$$

for suitable $A, B \in S$ and $\alpha, \beta \in \Gamma$. Since $J$ is a two-sided $\Gamma$-ideal, each of these expressions belongs to $J$ whenever $X_k \in J$. Hence $X_{k+1} \in J$, completing the induction. ■

Chemically, this proposition states that once a system enters a two-sided $\Gamma$-ideal, all states reachable via reaction pathways remain confined within that ideal. Thus, two-sided $\Gamma$-ideals model reaction basins or domains in which the chemistry is dynamically trapped under the available mediators.

## 4.3 Prime and semiprime chemical ideals

To understand how reaction activity distributes across the state space, it is useful to introduce notions of primeness and semiprimeness that adapt classical ideal-theoretic concepts to the ternary Γ-setting.

**Definition 6** (Prime chemical ideal)**.** A proper chemical ideal $P \subsetneq S$ is called *prime* if whenever

$$[A, \alpha, B, \beta, C] \in P$$

for some $A, B, C \in S$ and $\alpha, \beta \in \Gamma$, then at least one of $A, B, C$ lies in $P$.

This notion extends the classical understanding of prime ideals in semiring theory, where primeness forbids internal factorization of elements outside the ideal (compare [4]).

Thus, a prime chemical ideal cannot contain the result of a mediated interaction without "detecting" the presence of one of its participants. Chemically, $P$ behaves like a region of the state space whose boundary is sufficiently sharp that it cannot be entered as the product of a reaction between three states all lying outside $P$. In this sense, a prime chemical ideal captures a subsystem in which one interaction or family of interactions dominates access to its interior.

**Definition 7** (Semiprime chemical ideal)**.** A chemical ideal $I \subseteq S$ is called *semiprime* if for every $A \in S$ the following implication holds: if

$$[A, \alpha, A, \beta, A] \in I \quad \text{for all } \alpha, \beta \in \Gamma,$$

then $A \in I$.

This condition generalizes the classical semiprime property in semiring theory, where self-combinations inside an ideal force membership of the element itself (see [4]).

Here, $[A, \alpha, A, \beta, A]$ may be understood as a self-interaction or self-combination of the state $A$ under all possible mediators. The definition says that if every such self-interaction of $A$ falls inside $I$, then $A$ itself must already belong to $I$. Semiprimeness thus prevents the existence of "hidden" states outside $I$ whose entire mediated self-dynamics is trapped within $I$.

**Proposition 3.** *Every prime chemical ideal is semiprime.*

*Proof.* Let $P$ be a prime chemical ideal and suppose that $[A, \alpha, A, \beta, A] \in P$ for all $\alpha, \beta \in \Gamma$. In particular, there exist $\alpha_0, \beta_0 \in \Gamma$ such that $[A, \alpha_0, A, \beta_0, A] \in P$. By primeness, at least one of the three entries in this interaction must belong to $P$. Since all three are equal to $A$, we conclude that $A \in P$. Thus $P$ is semiprime. ■

From a chemical perspective, this result indicates that in a prime chemical ideal, any state whose self-interactions are entirely trapped within the ideal must itself be regarded as belonging to that ideal. Prime subsystems therefore exclude the possibility of persistent external states whose internal dynamics is indistinguishable, in terms of reaction products, from that of genuine internal states.

The theory of prime and semiprime chemical ideals provides a way to decompose a TGS-chemical system into structurally meaningful components, reflecting how reaction activity and mediated transformations are distributed across the state space. Comparable decomposition principles appear in algebraic treatments of semirings and their ideal lattices [6]. . In subsequent work, one may associate to a given system an appropriate spectrum of prime chemical ideals and study its topology, thus connecting the present framework with geometric methods.

# 5 Homomorphisms of chemical TGS

Homomorphisms provide a natural mechanism for comparing different TGS-chemical systems and transporting reaction behavior from one system to another. Just as homomorphisms of semirings or semigroups preserve algebraic structure, homomorphisms of ternary $\Gamma$-semirings preserve the mediated ternary interaction that encodes chemical transformation. The present section formalizes this notion and explains its chemical significance.

Throughout, $(S, \Gamma, [\,])$ and $(S', \Gamma, [\,]')$ denote two TGS-chemical systems sharing the same parameter set $\Gamma$. The requirement of a common $\Gamma$ reflects that mediators (catalysts, solvents, environmental conditions) are interpreted as parameters intrinsic to the interaction law and therefore must be preserved.

## 5.1 Definition

**Definition 8.** A map

$$f\colon S \longrightarrow S'$$

is called a *homomorphism of TGS-chemical systems* if for all $A, B, C \in S$ and all $\alpha, \beta \in \Gamma$,

$$f([A, \alpha, B, \beta, C]) \;=\; [f(A), \alpha, f(B), \beta, f(C)]'.$$

This condition is analogous to structure-preserving maps in classical semiring and algebraic systems, where homomorphisms preserve the underlying interaction laws (see [4, 6]).

Thus, $f$ commutes with the ternary $\Gamma$-operation: applying the reaction operation in $S$ and then mapping the result via $f$ yields the same state as first mapping the inputs via $f$ and then applying the reaction operation in $S'$. In other words, $f$ is a structure-preserving transformation of chemical environments.

Several immediate properties follow directly from the definition.

**Proposition 4.** *Let $f\colon S \to S'$ be a TGS-homomorphism.*

(a) *If $R \subseteq S$ is reaction-closed, then $f(R)$ is reaction-closed in $S'$.*

(b) *If $I \subseteq S$ is a chemical ideal, then $f(I)$ is a chemical ideal in $S'$.*

(c) *If $J \subseteq S$ is a $\Gamma$-ideal of any type (left, right, middle, or two-sided), then $f(J)$ is a $\Gamma$-ideal of the corresponding type in $S'$.*

*Proof.* Each property is verified by direct substitution using the homomorphism identity. For example, if $A, B, C \in R$ and $\alpha, \beta \in \Gamma$, then the reaction-closedness of $R$ gives $[A, \alpha, B, \beta, C] \in R$, and applying $f$ yields

$$f([A, \alpha, B, \beta, C]) = [f(A), \alpha, f(B), \beta, f(C)]' \in f(R),$$

establishing reaction-closedness of $f(R)$. The remaining cases follow the same pattern. ■

This result shows that homomorphisms are compatible with the structural subsystems developed in Section 4: reactors, basins, and pathways are mapped to reactors, basins, and pathways in the target system.

## 5.2 Chemical meaning

A TGS-homomorphism models a consistency-preserving transformation between chemical environments. Its chemical interpretations include the following:

**(1) Change of solvent or medium.** Suppose $S$ describes reaction behavior in solvent $X$ and $S'$ in solvent $Y$. A homomorphism $f\colon S \to S'$ represents a map translating chemical states from the $X$-environment to the $Y$-environment such that the mediated interactions are preserved: a triple interaction in $X$ corresponds exactly to the mapped triple interaction in $Y$. This formalizes the intuitive idea that a well-defined solvent change should send reaction pathways to reaction pathways without altering their essential structure.

**(2) Change of catalyst or catalytic regime.** Different catalytic environments can be modeled by different TGS-chemical systems built on the same parameter space $\Gamma$ but with distinct state spaces or distinct ternary interaction laws. A homomorphism

$$f\colon S \to S'$$

can represent the adjustment of reaction behavior when switching from one catalyst to another. The preservation of the mediated operation ensures that catalytic effects are transferred systematically rather than arbitrarily.Such environment-to-environment mappings have analogues in algebraic treatments of ternary and parameter-dependent transformations (see [13]).

**(3) Controlled mapping between chemical environments.** More generally, a homomorphism encodes any structured change of environment where reaction behavior is transformed coherently. This may represent, for example:

- embedding a system with restricted state space into a larger one;
- coarse-graining a complex reaction network into a simpler model;
- mapping between different thermodynamic or field-controlled environments;
- abstraction from microscopic to effective macroscopic states.

In each of these examples, the homomorphism ensures that reaction mechanisms and mediator influences retain their form under translation.

**(4) Compatibility with reaction pathways.** Since homomorphisms preserve the ternary Γ-operation, they also preserve reaction pathways in the sense of Section 4. Every reaction sequence

$$X_0 \to X_1 \to \cdots \to X_n$$

in $S$ is carried by $f$ to a reaction pathway

$$f(X_0) \to f(X_1) \to \cdots \to f(X_n)$$

in $S'$. Thus, homomorphisms provide a bridge between dynamical behaviors in different systems, enabling the systematic study of how pathways transform under environmental changes.

Overall, homomorphisms of TGS-chemical systems play a role analogous to structure-preserving maps in algebra, but their chemical interpretation is richer: they express how reaction laws, mediators, and transformation dynamics behave under coherent changes of environment. This makes them powerful tools for both mathematical analysis and chemical modeling.

# 6 Examples from chemistry

In this section we present several abstract but chemically meaningful examples illustrating how mediated ternary interactions arise naturally in chemical systems.In this section we present several abstract but chemically meaningful examples illustrating how mediated ternary interactions arise naturally in chemical systems. These examples parallel mathematically formal approaches to chemical structure and transformation found in classical mathematical chemistry (see [1, 8, 12]).

The purpose of these examples is not to model specific experimental systems but to show how familiar chemical phenomena can be expressed within the TGS framework introduced above.

Throughout, $(S, \Gamma, [\,])$ denotes a TGS-chemical system in the sense of Section 3, where $S$ represents chemical states and $\Gamma$ represents mediating conditions.

## 6.1 Catalyzed reactions

Catalysis provides a direct example of a mediated transformation in which the presence of a catalyst modifies the reaction pathway without being consumed. Let $A, B, C \in S$ denote chemical states that participate in a multi-step reaction, and let $\alpha, \beta \in \Gamma$ represent catalytic regimes.

Consider the ternary operation

$$[A, \alpha, B, \beta, C] = D.$$

Here, $A$ may be interpreted as an initial reactant state, $B$ as an interacting partner or intermediate, and $C$ as a subsequent state through which the system passes. The mediator $\alpha$ can encode the presence of a catalyst that opens a specific reaction pathway, while $\beta$ may represent a co-catalyst or a secondary catalytic condition.

If $\alpha$ corresponds to a catalyst that lowers the effective barrier between $A$ and $B$, and $\beta$ indexes a catalytic effect acting on the transformation from $B$ to $C$, then $D$ represents the resulting state of the catalyzed sequence.

Different choices of $\alpha$ and $\beta$ generally produce different outcomes:

$$[A, \alpha_1, B, \beta, C] \neq [A, \alpha_2, B, \beta, C],$$

even when the underlying species $A, B, C$ remain fixed. This expresses, in algebraic form, the well-known fact that changing catalysts can modify the reaction pathway or final products while preserving stoichiometry.

The ternary structure is essential here: the catalyst is not appended externally but serves as an intrinsic argument of the reaction law.

## 6.2 Phase transitions under thermodynamic control

Phase transformations depend sensitively on thermodynamic parameters such as temperature and pressure. In the TGS framework, such environmental conditions are naturally represented as elements of $\Gamma$.

Let $\Gamma$ consist of pairs $(T, p)$ corresponding to permissible temperature–pressure regimes.Let $A, B, C \in S$ represent physical states of a substance, such as configurations or phase descriptors. A ternary interaction

$$[A, (T_1, p_1), B, (T_2, p_2), C]$$

produces a state $D$, where the mediators $(T_1, p_1)$ and $(T_2, p_2)$ govern the transitions between $A \to B$ and $B \to C$ respectively.

For example:

- If $(T_1, p_1)$ represents conditions favoring melting, and $(T_2, p_2)$ represents conditions favoring vaporization, then $D$ may correspond to a higher-energy phase.

- If $(T_1, p_1)$ lies in a stability region for a solid phase, and $(T_2, p_2)$ lies in a stability region for a metastable phase, then $D$ may encode a metastable state reached by sequential transitions.

The ternary formulation captures the fact that multi-step phase transformations are governed not only by initial and final conditions but also by intermediate thermodynamic regimes. Different paths through $(T, p)$- space yield different outcomes, and the dependence is faithfully recorded by the mediators in the Γ-operation.

## 6.3 Quantum state transitions under external fields

Quantum systems subject to external electromagnetic fields provide a further setting in which ternary, parameter-dependent interactions arise. Let $S$ denote a set of quantum states, which may include electronic, vibrational, or spin configurations. Let $\Gamma$ index external field parameters such as field strength, frequency, or polarization.Such field-mediated transitions have abstract algebraic analogues in parameter-dependent ternary relation frameworks (see [13]).

A ternary interaction

$$[A, \alpha, B, \beta, C] = D$$

may then model a sequence of field-induced transitions:

- $A \to B$ mediated by field parameter $\alpha$,
- $B \to C$ mediated by field parameter $\beta$,
- resulting in a state $D$ after the composite process.

For instance:

- $\alpha$ may represent a low-frequency field inducing a transition from $A$ to $B$;
- $\beta$ may represent a high-frequency field inducing a transition from $B$ to $C$;
- the final state $D$ depends on the combined effect of both fields in sequence.

The value of $D$ may differ significantly from what is obtained by either field alone, reflecting the well-established sensitivity of quantum transitions to external field combinations. The ternary structure captures this dependence by integrating the field parameters directly into the reaction law.

These examples illustrate how the ternary Γ-operation provides a natural formalism for expressing catalysis, thermodynamic control, and field-induced quantum transitions within a single coherent algebraic framework. The examples are intentionally abstract, focusing on the structural features that make TGS-chemical systems flexible enough to encode a wide range of chemical behavior.

# 7 Conclusion and future work

In this work we have developed an axiomatic framework for modeling chemical systems using the structure of a ternary Γ-semiring. Beginning with the observation that chemical transformations are inherently multi-parameter and multi-state processes, we formulated a reaction law in which chemical states and mediating conditions appear as intrinsic arguments of a ternary operation. This contrasts with the classical binary perspective, where catalysts and environmental factors are appended externally rather than participating structurally in the transformation process.

The foundational contribution of the paper lies in isolating the mathematical axioms that govern such mediated transformations and demonstrating how these axioms admit chemically meaningful interpretations. These axioms extend the structural principles familiar from classical semiring theory (see [3,4]) to a ternary Γ-mediated setting appropriate for chemical applications.

The ternary operation encodes multi-step transformations, the Γ-parameters incorporate catalytic and environmental effects, and the associativity and distributivity relations reflect coherence of reaction pathways. The resulting concept of a TGS-chemical system provides a unified formalism in which multi-state, catalyst-dependent, and environment-dependent phenomena can be described algebraically.

We also developed the structural theory of these systems, introducing chemical ideals, Γ-ideals, and their prime and semiprime variants. These notions identify chemically stable subsystems, reaction-closed domains, and regions whose internal behavior governs the structure of mediated interactions. Reaction pathways were characterized in terms of iterated ternary operations, and we showed how homomorphisms between TGS-chemical systems provide consistency-preserving maps between different chemical environments. Finally, we illustrated the framework with abstract examples drawn from catalysis, phase transitions, and field-driven quantum processes.

## Future directions

The present framework opens several avenues for further development.

- **Kinetic and dynamical refinements.** While the TGS formalism captures structural relationships between states and mediators, incorporating explicit temporal or kinetic data would allow the construction of mediated dynamical systems. A natural direction is to study sequences of ternary interactions as discrete dynamical processes and to identify stability, periodicity, or convergence phenomena within this setting.

- **Quantitative extensions.** The current theory treats $S$ and $\Gamma$ abstractly. Enriching these sets with additional algebraic or topological structure—such as orders, metrics, or weights—could allow the encoding of reaction energetics, field strengths, or graded catalytic effects. Such extensions would be essential for connecting the TGS framework to numerical models.

- **Categorical and geometric viewpoints.** The ideal theory developed here suggests the possibility of defining spectra of prime chemical ideals and studying their geometric features. This may lead to a form of *ternary $\Gamma$-geometry* in which chemical structure is represented through geometric invariants of the spectrum.

- **Computational and AI-based models.** The unified treatment of states and mediators makes TGS-chemical systems natural candidates for symbolic or rule-based computational models. Subsequent work may explore how the ternary operation interacts with algorithmic reasoning, abstraction, or machine-assisted simulation, thereby linking algebraic chemistry with emerging methodologies in computational chemistry and symbolic AI.Such directions resonate with semiring-based computational frameworks and ternary parameterized transformations explored in abstract algebraic settings (see [13]).

Overall, the ternary Γ-semiring viewpoint offers a flexible and conceptually coherent foundation for the algebraic study of chemical systems. The theory developed in this paper establishes the basic structure on which further analytical, geometric, and computational developments can be built.

***Acknowledgment*:** Chandrasekhar Gokavarapu gratefully acknowledges the supervision of Dr. D. Madhusudhana Rao in the conceptualization and development of this work. The authors express their sincere thanks to Dr. Ramachandra R. K., Principal, Government College (Autonomous), Rajahmundry, for providing continuous support and an excellent academic environment. The authors also acknowledge the institutional support of Acharya Nagarjuna University, Guntur. They further thank their colleagues in the Departments of Mathematics and Chemistry for helpful discussions during the preparation of this manuscript.